\documentclass[12pt]{amsart}
\usepackage{amsmath, amsthm, amssymb}
\title{On Uniform Definability of Types over Finite Sets}
\author{Vincent Guingona* \\ Department of Mathematics \\ University of Maryland }
\thanks{*Partially Supported by Laskowski's NSF grants DMS-0600217 and 0901336}
\date{\today}

\newtheorem{thm}{Theorem}[section]
\newtheorem{cor}[thm]{Corollary}
\newtheorem{lem}[thm]{Lemma}
\newtheorem{prop}[thm]{Proposition}

\newtheorem{ques}[thm]{Open Question}

\theoremstyle{remark}
\newtheorem{rem}[thm]{Remark}

\theoremstyle{definition}
\newtheorem{defn}[thm]{Definition}

\newcommand{\dom}{\mathrm{dom} }
\newcommand{\tp}{\mathrm{tp} }
\newcommand{\concat}{{}^{\frown} }
\newcommand{\ID}{\mathrm{ID} }

\newcommand{\leng}{\mathrm{lg} }
\newcommand{\opp}{\mathrm{opp} }

\newcommand{\Th}{\mathrm{Th} }

\begin{document}

\begin{abstract}
 In this paper, using definability of types over indiscernible sequences as a template, we study a property of formulas and theories called ``uniform definability of types over finite sets'' (UDTFS).  We explore UDTFS and show how it relates to well-known properties in model theory.  We recall that stable theories and weakly o-minimal theories have UDTFS and UDTFS implies dependence.  We then show that all dp-minimal theories have UDTFS.
\end{abstract}

\maketitle

\section{Introduction}\label{Section_intro}

The notion of definability of types is a useful tool in the study of model theory.  The nature of definability of $\varphi$-types for stable formulas $\varphi(\overline{x}; \overline{y})$ is well understood, but we wish to extend this notion to the more general context of dependent formulas.  Restricting our attention to types over indiscernible sequences has a tendency to smooth things out.  We show that the existence of a uniform definition for $\varphi$-types over arbitrary indiscernible sequences characterizes stability and the existence of a uniform definition for $\varphi$-types over finite indiscernible sequences characterizes dependence.  The main question follows from this concept: If we replace ``indiscernible sequences'' with general sets, do these characterizations still hold?

The characterization of stable formulas in terms of definability of types over arbitrary sets certainly holds.  A formula $\varphi(\overline{x}; \overline{y})$ is stable if and only if there exists some other formula $\psi(\overline{y}; \overline{z})$ such that, for all $\varphi$-types $p(\overline{x})$, there exists a tuple $\overline{c}$ from $\dom(p)$ so that $\psi(\overline{y}; \overline{c})$ defines the type $p$ (see \cite{shelah}).  Using this and the indiscernible sequence analogy as a template, we define a new notion of definability of types called ``uniform definability of types over finite sets'' (UDTFS).  Instead of considering all $\varphi$-types, $p$, we consider only those $\varphi$-types that have finite domain.  With this one simple modification, we expand the class of theories that have definability of types well beyond stable theories.  In fact, it is conjectured that this notion actually characterizes dependence.

This paper exhibits some previously known results about UDTFS including basic properties and its relation to some other model-theoretic dividing lines.  In particular, we show that stable theories and weakly o-minimal theories have UDTFS and that UDTFS theories are dependent.  These relations have been worked out by Hunter Johnson and Chris Laskowski \cite{LJ2009}, but are included here for completeness.  We then proceed to prove a new generalization of the fact that weakly o-minimal theories have UDTFS; specifically, we show that dp-minimal theories have UDTFS.

We would like to thank several people for helping to make this paper possible.  We are especially grateful to Chris Laskowski for his insights into the stable case of Theorem \ref{Thm_IndiscernibleSequenceChar} and for suggesting ways to make the paper more presentable and understandable.  We also thank Alfred Dolich, John Goodrick, and Sergei Starchenko for listening to some preliminary versions of the arguments in this paper and giving helpful feedback.

\subsection{Definitions}

For this paper, a ``formula'' means a $\emptyset$-definable formula in a fixed language $L$ unless otherwise specified.  If $\theta(\overline{x})$ is a formula, then denote $\theta(\overline{x})^0 = \neg \theta(\overline{x})$ and $\theta(\overline{x})^1 = \theta(\overline{x})$.  Most formulas we work with are partitioned formulas, $\varphi(\overline{x}; \overline{y})$, where the variables are broken into two distinct sets.  The first set of variables is called the set of \textit{free variables} and the second set is called the set of \textit{parameter variables}.  When we have a list of tuples of variables, we will sometimes denote this with a boldface variable to shorten notation.  For example, we could write $\varphi(\overline{x}; \overline{y}_0, ..., \overline{y}_{n-1})$ as $\varphi(\overline{x}; \overline{\mathbf{y}})$.

We work in a complete first-order theory $T$ in the language $L$ with a monster model, $\mathfrak{C}$.  Fix a partitioned $L$-formula $\varphi(\overline{x}; \overline{y})$.  We say that a set $B \subseteq \mathfrak{C}^{\leng(\overline{y})}$ is $\varphi$-\textit{independent} if, for all maps $s \in {}^B 2$, the set of formulas $\{ \varphi(\overline{x}; \overline{b})^{s(\overline{b})} : \overline{b} \in B \}$ is consistent.  We say that $\varphi$ has \textit{independence dimension} $N < \omega$, which we denote by $\ID(\varphi) = N$, if $N$ is maximal such that there exists a $\varphi$-independent set $B \subseteq \mathfrak{C}^{\leng(\overline{y})}$ with $|B| = N$.  If no such maximal $N$ exists, we say that $\ID(\varphi) = \infty$ and $\varphi$ is \textit{independent}.  We say that $\varphi$ is \textit{dependent} (some authors call this NIP for ``not the independence property'') if $\ID(\varphi) = N$ for some $N < \omega$.  Finally, we say that a theory $T$ is \textit{dependent} if all partitioned formulas are dependent.

By a ``$\varphi$-type over $B$'' for some small subset $B \subseteq \mathfrak{C}^{\leng(\overline{y})}$ we mean a maximal consistent set of formulas of the form $\varphi(\overline{x}; \overline{b})^t$ for some $t < 2$ and some $\overline{b} \in B$.  If $p$ is a $\varphi$-type over $B$, then we say that $p$ has \textit{domain} $\dom(p) = B$.  For any small subset $B \subseteq \mathfrak{C}^{\leng(\overline{y})}$, the space of all $\varphi$-types with domain $B$ is denoted $S_\varphi(B)$.  Any $\varphi$-type $p$ with domain $B$ gives rise to a function $\delta : B \rightarrow 2$ where, for all $\overline{b} \in B$, $\varphi(\overline{x}; \overline{b})^{\delta(\overline{b})} \in p(\overline{x})$.  Therefore,
\[
 p(\overline{x}) = \{ \varphi(\overline{x}; \overline{b})^{\delta(\overline{b})} : \overline{b} \in B \}.
\]
We will call this $\delta$ the function \textit{associated} to the $\varphi$-type $p$.  For $B_0 \subseteq B$, $p \in S_\varphi(B)$, and $\delta$ associated to $p$, let
\[
 p_{B_0} (\overline{x}) = \{ \varphi(\overline{x}; \overline{b})^{\delta(\overline{b})} : \overline{b} \in B_0 \}
\]
denote the restriction of $p$ to $B_0$.  For sets $B_1 \subseteq B_0 \subseteq B$, let
\[
 p_{B_0,B_1} (\overline{x}) = \{ \varphi(\overline{x}; \overline{b})^{\delta(\overline{b})} : \overline{b} \in B_0 - B_1 \} \cup \{ \neg \varphi(\overline{x}; \overline{b})^{\delta(\overline{b})} : \overline{b} \in B_1 \}.
\]
That is, $p_{B_0,B_1}$ is $p_{B_0}$ except we negate all instances of $\varphi$ on elements of $B_1$.  Sometimes we call this perturbing $p_{B_0}$ by $B_1$.  One should note that $p_{B_0,B_1}$ need not be a $\varphi$-type because it need not be consistent.

Consider the following notion of definability of types:

\begin{defn}\label{Def_Defines}
 Fix a $\varphi$-type $p(\overline{x})$ and a formula $\psi(\overline{y})$ that is not necessarily over $\emptyset$.  We say that $\psi$ \textit{defines} $p$ if, for all $\overline{b} \in \dom(p)$, $\models \psi(\overline{b})$ if and only if $\varphi(\overline{x}; \overline{b}) \in p(\overline{x})$.
\end{defn}

One should note that it is easy, in general, to find a definition of a $\varphi$-type, $p$.  For example, if $\overline{a}$ is any realization of $p(\overline{x})$, then we can simply take $\psi(\overline{y}) = \varphi(\overline{a}; \overline{y})$ as a definition of $p$.  The difficulty comes in finding a definition that is defined over $\dom(p)$.  In fact, the existence of a defining formula over $\dom(p)$ for all $\varphi$-types $p$ actually characterizes the stability of $\varphi$ (see \cite{shelah}).  To motivate the remaining sections of this paper, we now consider the characterizations of stability and dependence in terms of the existence of uniform definitions for $\varphi$-types over indiscernible sequences.

\subsection{Indiscernible Sequences and Definability of Types}

When considering indiscernible sequences in relation to a partitioned formula $\varphi(\overline{x}; \overline{y})$, it helps to look at the following set of formulas:

\begin{equation}\label{Eq_defofDeltaN}
 \Delta_N(\overline{y}_0, ..., \overline{y}_{N-1}) = \left\{ \exists \overline{x} \bigwedge_{i < N} \varphi(\overline{x}; \overline{y}_i)^{s(i)} : s \in {}^N 2 \right\}.
\end{equation}

More generally, fix any collection of formulas $\Delta(\overline{y}_0, ..., \overline{y}_N)$ and $(I, <)$ a linear order.  A $\Delta$-\textit{indiscernible sequence} is a sequence $\langle \overline{b}_i : i \in I \rangle$ such that, for all $i_0 < ... < i_N$ and $j_0 < ... < j_N$ from $I$, for all $\delta \in \Delta$, $\models \delta(\overline{b}_{i_0}, ..., \overline{b}_{i_N})$ if and only if $\models \delta(\overline{b}_{j_0}, ..., \overline{b}_{j_N})$.  We say that a $\Delta$-indiscernible sequence $\langle \overline{b}_i : i \in I \rangle$ is a $\Delta$-\textit{indiscernible set} if, for all $i_0, ..., i_N \in I$ distinct and all $j_0, ..., j_N \in I$ distinct (without any restriction on the ordering), for all $\delta \in \Delta$, $\models \delta(\overline{b}_{i_0}, ..., \overline{b}_{i_N})$ if and only if $\models \delta(\overline{b}_{j_0}, ..., \overline{b}_{j_N})$.  We simply say that $\langle \overline{b}_i : i \in I \rangle$ is an indiscernible sequence (respectively, indiscernible set) if it is a $\Delta$-indiscernible sequence (respectively, $\Delta$-indiscernible set) for all sets of formulas, $\Delta$ (with a partitioned of the free variables into $\leng(\overline{b}_i)$-tuples).  For any set of formulas $\Delta$, let $\pm \Delta = \{ \delta^t : \delta \in \Delta, t < 2 \}$.  We say that a set of formulas $\Delta(\overline{y}_0, ..., \overline{y}_N)$ is \textit{closed under permutations} if, for all $\delta \in \Delta$ and all $\sigma \in S_{N+1}$ (where $S_n$ is the group of permutations on $n$), the formula
\[
 \delta(\overline{y}_{\sigma(0)}, ..., \overline{y}_{\sigma(N)})
\]
is also in $\Delta$.  For example, $\Delta_N$ as in \eqref{Eq_defofDeltaN} is closed under permutations.

\begin{thm}\label{Thm_IndiscernibleSequenceChar}
 For a partitioned formula, $\varphi(\overline{x}; \overline{y})$, the following hold:

 \begin{itemize}
  \item [(i)] $\varphi$ is stable if and only if there exists $\psi(\overline{y}; \overline{z}_0, ..., \overline{z}_{K-1})$ such that, for all indiscernible sequences $\langle \overline{b}_i : i \in I \rangle$ with $|I| \ge 2$ and all $p(\overline{x}) \in S_\varphi(\{ \overline{b}_i : i \in I \})$, there exists $i_0, ..., i_{K-1} \in I$ such that $\psi(\overline{y}; \overline{b}_{i_0}, ..., \overline{b}_{i_{K-1}})$ defines $p(\overline{x})$.
  \item [(ii)] $\varphi$ is dependent if and only if there exists $\psi(\overline{y}; \overline{z}_0, ..., \overline{z}_{K-1})$ such that, for all finite indiscernible sequences $\langle \overline{b}_i : i \in L \rangle$ with $L \ge 2$ and all $p(\overline{x}) \in S_\varphi(\{ \overline{b}_i : i \in L \})$, there exists $i_0, ..., i_{K-1} \in L$ such that $\psi(\overline{y}; \overline{b}_{i_0}, ..., \overline{b}_{i_{K-1}})$ defines $p(\overline{x})$.
 \end{itemize}
\end{thm}

Before giving the proof of Theorem \ref{Thm_IndiscernibleSequenceChar}, let us first deal with the case where $\varphi$ is dependent.  If $\Delta$ is any set of formulas that is closed under permutations, we show that, for any $\Delta$-indiscernible sequence that is not a $\Delta$-indiscernible set, there exists an instance of $\pm \Delta$ that defines the linear order of the indiscernible sequence.  The proof of this lemma is based on the proof of Theorem II.4.7 in \cite{shelah}.

\begin{lem}\label{ordersensitiveformula}
 If $\Delta(\overline{y}_0, ..., \overline{y}_N)$ is a set of formulas that is closed under permutations, $(I, <)$ is a linear order, and $\langle \overline{b}_i : i \in I \rangle$ is a $\Delta$-indiscernible sequence that is not a $\Delta$-indiscernible set, then there exists $t < n - 1$ and $\delta \in \pm \Delta$ so that
 \[
  \models \delta(\overline{b}_{i_0}, ..., \overline{b}_{i_{t-1}}, \overline{b}_{i_t}, \overline{b}_{i_{t+1}}, ..., \overline{b}_{i_N}) \wedge \neg \delta(\overline{b}_{i_0}, ..., \overline{b}_{i_{t-1}}, \overline{b}_{i_{t+1}}, \overline{b}_{i_t}, \overline{b}_{i_{t+2}}, ..., \overline{b}_{i_N}).
 \]
 for some (equivalently all) $i_0 < i_1 < ... < i_N$ from $I$.  That is, $\delta$ is ``order sensitive'' at $t$.
\end{lem}

\begin{proof}
 To simplify notation, assume that $0 < 1 < ... < N$ is in $I$ and show this for $\overline{b}_0, ..., \overline{b}_N$.  Since $\langle \overline{b}_i : i \in I \rangle$ is a $\Delta$-indiscernible sequence that is not a $\Delta$-indiscernible set, there exists some $\delta' \in \pm \Delta$ witnessing this fact.  That is, $\models \delta'(\overline{b}_0, ..., \overline{b}_N)$ but, for some $\sigma \in S_{N+1}$, $\models \neg \delta'(\overline{b}_{\sigma(0)}, ..., \overline{b}_{\sigma(N)})$.  Now $S_{N+1}$, as a group, is generated by permutations of the form $(t \ t\text{+1})$ for $t < N$ (the permutation that is the identity on all of $N+1$ except it swaps $t$ and $t+1$).  Therefore, there exists $\sigma' \in S_{N+1}$ and $t < N$ such that
 \[
  \models \delta'(\overline{b}_{\sigma'(0)}, ..., \overline{b}_{\sigma'(N)}) \wedge \neg \delta'(\overline{b}_{\tau \circ \sigma'(0)}, ..., \overline{b}_{\tau \circ \sigma'(N)})
 \]
 where $\tau = (t \ t\text{+1})$.  Since $\Delta$ is closed under permutations, if we let
 \[
  \delta(\overline{y}_0, ..., \overline{y}_{N}) = \delta'(\overline{y}_{\sigma'(0)}, ..., \overline{y}_{\sigma'(N)}),
 \]
 then $\delta \in \pm \Delta$.  Therefore, we get
 \[
  \models \delta(\overline{b}_0, ..., \overline{b}_{t-1}, \overline{b}_t, \overline{b}_{t+1}, ..., \overline{b}_{N}) \wedge \neg \delta(\overline{b}_0, ..., \overline{b}_{t-1}, \overline{b}_{t+1}, \overline{b}_t, \overline{b}_{t+2}, ..., \overline{b}_{N})
 \]
 as desired.
\end{proof}

If $(I, <) = (L, <)$ is finite but $L > N$, then we can take the initial $t$ elements and final $N-t-1$ elements of $\langle \overline{b}_i : i \in L \rangle$ and we get that the formula 
\[
 \theta(\overline{y}_0; \overline{y}_1) = \delta(\overline{b}_0, ..., \overline{b}_{t-1}, \overline{y}_0, \overline{y}_1, \overline{b}_{L-N+t+1}, ..., \overline{b}_{L-1})
\]
defines the linear order of the indiscernible sequence on $\langle \overline{b}_i : t \le i \le L-N+t \rangle$.  That is, for all distinct $i, j$ with $t \le i,j \le L-N+t$, $\models \theta(\overline{b}_i; \overline{b}_j)$ if and only if $i < j$.  We use this definition of the linear order, in conjunction with the following lemma, to get a definition of $\varphi$-types over finite indiscernible sequences for dependent $\varphi$.  The proof of the following lemma can be found in \cite{shelah} (see Theorem II.4.13).

\begin{lem}\label{Lemma_Alternation}
 Fix a dependent partitioned formula $\varphi(\overline{x}; \overline{y})$, let $N = \ID(\varphi)$, and let $\Delta = \Delta_{N+1}$ as in \eqref{Eq_defofDeltaN}.  If $(I, <)$ is a linear order, $\langle \overline{b}_i : i \in I \rangle$ is a $\Delta$-indiscernible sequence, and $\overline{a} \in \mathfrak{C}^{\leng(\overline{x})}$, then there exists $K \le N+1$ $<$-convex subsets of $I$, $I_0, ..., I_{K-1}$, such that $\models \varphi(\overline{a}; \overline{b}_i)$ if and only if $i \in I_0 \cup ... \cup I_{K-1}$.  If $\langle \overline{b}_i : i \in I \rangle$ is a $\Delta$-indiscernible set, then either $|I| \le 2N+1$ or there exists $K \le N$ elements $i_0, ..., i_{K-1} \in I$ and $t < 2$ such that $\models \varphi(\overline{a}; \overline{b}_i)^t$ if and only if $i \in \{ i_0, ..., i_{K-1} \}$.
\end{lem}

We are now ready to prove the main motivating theorem.  Fix any partitioned formula, $\varphi(\overline{x}; \overline{y})$.

\begin{proof}[Proof of Theorem \ref{Thm_IndiscernibleSequenceChar}]
 (i): Suppose $\varphi$ is unstable.  Since $\varphi(\overline{x}; \overline{y})$ is unstable, it has the order property.  By compactness and Ramsey's Theorem, there exists an indiscernible sequence $\langle \overline{b}_q : q \in \mathbb{Q} \rangle$ and a set $\{ \overline{a}_r : r \in \mathbb{R} \}$ such that $\models \varphi(\overline{a}_r; \overline{b}_q)$ if and only if $r < q$.  If there existed a formula $\psi(\overline{y}; \overline{z}_0, ..., \overline{z}_{K-1})$ that uniformly defined all $\varphi$-types over indiscernible sequences, then each $\overline{a}_r$ would be determined by an element in $( \{ \overline{b}_q : q \in \mathbb{Q} \} )^K$, which is countable, hence $| \{ \overline{a}_r : r \in \mathbb{R} \} | = \aleph_0$.  However, each $\overline{a}_r$ is distinct, hence $| \{ \overline{a}_r : r \in \mathbb{R} \} | = 2^{\aleph_0}$.  Therefore, no such formula $\psi$ exists.

 Conversely, if we assume that $\varphi$ is stable, then Shelah showed the existence of a formula $\psi(\overline{y}; \overline{z}_0, ..., \overline{z}_{K-1})$ that uniformly defines $\varphi$-types over all sets in Theorem II.2.12 of \cite{shelah}, so $\psi$ certainly works over indiscernible sequences.

 (ii): Assume $\varphi$ is independent.  By compactness and Ramsey's theorem, there exists $\langle \overline{b}_i : i < \omega \rangle$ that is both indiscernible and $\varphi$-independent.  If there existed a $\psi(\overline{y}; \overline{z}_0, ..., \overline{z}_{K-1})$ that uniformly defines all $\varphi$-types over finite indiscernible sequences, then for each $L < \omega$, the number of $\varphi$-types over $\{ \overline{b}_i : i < L \}$ would be bounded by $L^K$.  However, since $\{ \overline{b}_i : i < L \}$ is $\varphi$-independent, the number of $\varphi$-types over $\{ \overline{b}_i : i < L \}$ is exactly $2^L$.  This gives us that $2^L \le L^K$ for all $L < \omega$.  However $K$ is fixed, so the exponential function will eventually overtake the polynomial function, giving us a contradiction.

 Conversely, suppose $\varphi$ is dependent.  Let $N = \ID(\varphi)$ and let $\Delta = \Delta_{N+1}$ as in \eqref{Eq_defofDeltaN} (note that $\Delta$ is a finite set, hence so is $\pm \Delta$).  We define finitely many formulas, $\psi_\ell(\overline{y}; \overline{\mathbf{z}})$, so that, for all finite indiscernible sequences $\langle \overline{b}_i : i \in L \rangle$ and all $p(\overline{x}) \in S_\varphi( \{ \overline{b}_i : i \in L \} )$, there exists $\ell$ and $i_0, ..., i_{K} \in L$ such that $\psi_\ell(\overline{y}; \overline{b}_{i_0}, ..., \overline{b}_{i_{K}})$ defines $p(\overline{x})$ (for appropriate $K$).  As in Lemma \ref{Lemma_ManyOneUDTFS} below, one can combine these formulas into a single $\psi(\overline{y}; \overline{\mathbf{z}})$, as desired.  First off, for the case where our indiscernible sequence $\langle \overline{b}_i : i \in L \rangle$ happens to be a $\Delta$-indiscernible set or where $L \le 2N+1$, we can use the defining formula
 \[
  \psi_t(\overline{y}; \overline{z}_0, ..., \overline{z}_{2N}) = \left( \bigvee_{i \le 2N} \overline{y} = \overline{z}_i \right)^t
 \]
 for some $t < 2$ by Lemma \ref{Lemma_Alternation}.  Otherwise, $\langle \overline{b}_i : i \in L \rangle$ is a $\Delta$-indiscernible sequence that is not a $\Delta$-indiscernible set.  So, by Lemma \ref{ordersensitiveformula}, there exists $t < N$ and $\delta \in \pm \Delta$ such that the formula
 \[
  \theta_{t,\delta}(\overline{y}_0, \overline{y}_1; \overline{z}_0, ..., \overline{z}_{t-1}, \overline{z}_{t+2}, ..., \overline{z}_N) = \delta(\overline{z}_0, ..., \overline{z}_{t-1}, \overline{y}_0, \overline{y}_1, \overline{z}_{t+2}, ..., \overline{z}_N)
 \]
 can be used to define the linear order of the indiscernible sequence.  By Lemma \ref{Lemma_Alternation}, for any $\overline{a}$, there exists $K \le N+1$ intervals of $I$, $I_0, ..., I_{K-1}$ so that $\models \varphi(\overline{a}; \overline{b}_i)$ if and only if $i \in I_0 \cup ... \cup I_{K-1}$.  So the truth value of $\varphi(\overline{a}; \overline{b}_i)$ for $t \le i \le L - N + t$ is determined by the formula 
 \[
  \theta_{t,\delta}(\overline{y}_0, \overline{y}_1; \overline{b}_0, ..., \overline{b}_{t-1}, \overline{b}_{L-N+t+1}, ..., \overline{b}_{L-1})
 \]
 used to define the ordering and at most $2N+2$ boundary points of the intervals, $I_0, ..., I_{K-1}$.  We now see that there are a bounded number of formulas, $\psi_\ell(\overline{y}; \overline{\mathbf{z}})$ which depend only on $t < N$, $\delta \in \pm \Delta$, $K < N+1$, and the various possible truth values of $\varphi(\overline{a}; \overline{b}_i)$ for $i < t$ and $i > L-N+t$, so that the $\psi_\ell(\overline{y}; \overline{\mathbf{z}})$ uniformly define $\varphi$-types over finite indiscernible sequences.
\end{proof}

We use Theorem \ref{Thm_IndiscernibleSequenceChar} as a template when formulating the notion of UDTFS.  In the stable setting, indiscernible sequences can be replaced with arbitrary sets and we still get a characterization of stability in terms of definability of types.  The desire is to also get a characterization of dependence by replacing the condition of ``finite indiscernible sequences'' with ``finite sets.''  Though it is still an open problem whether or not uniform definability of types over finite sets characterizes dependence, we still know that many subclasses of dependent theories have UDTFS.  We spend the remainder of this paper exploring how known dividing lines for theories relate to UDTFS.


\section{Basic properties of UDTFS}\label{Section_BasicUDTFS}

Let $T$ be a complete theory and let $\mathfrak{C}$ be a monster model for $T$.  As an analogy to the stable case and to Theorem \ref{Thm_IndiscernibleSequenceChar} above, consider the following definition:

\begin{defn}\label{Def_UDTFS}
 A partitioned formula $\varphi(\overline{x}; \overline{y})$ has \textit{uniform definablility of types over finite sets} (\textit{UDTFS}) if there exists a formula $\psi(\overline{y}; \overline{z}_0, ..., \overline{z}_{k-1})$ such that, for all finite sets of $\leng(\overline{y})$-tuples $B$ with $|B| \ge 2$ and for all $p(\overline{x}) \in S_{\varphi}(B)$, there exist $\overline{c}_0, ..., \overline{c}_{k-1} \in B$ such that $\psi(\overline{y}; \overline{c}_0, ..., \overline{c}_{k-1})$ defines $p(\overline{x})$.  We call such a $\psi$ a \textit{uniform definition of $\varphi$-types over finite sets}.  A theory has \textit{uniform definablility of types over finite sets} (\textit{UDTFS}) if every formula has UDTFS.
\end{defn}

First, note that UDTFS transfers between elementarily equivalent structures.  It is still open whether or not reducts of UDTFS theories have UDTFS; one could accidentally ``throw out'' the definition when taking the reduct.

Lemma \ref{Lemma_BooleanCombos} and Proposition \ref{Prop_UDTFSDep} below are due to Johnson and Laskowski and are proved in \cite{LJ2009}, but we include proofs of them here for completeness.

\begin{lem}\label{Lemma_BooleanCombos}
 The class of formulas that have UDTFS and have the same free variables $\overline{x}$ is closed under boolean combinations.
\end{lem}

\begin{proof}
 Fix $\varphi(\overline{x}; \overline{y})$ and $\psi(\overline{x}; \overline{z})$ and suppose these formulas have UDTFS.  Let
 \[
  \gamma_\varphi(\overline{y}; \overline{w}_0, ..., \overline{w}_{n-1}) \text{ and } \gamma_\psi(\overline{z}; \overline{v}_0, ..., \overline{v}_{n-1})
 \]
 be uniform definitions of $\varphi$-types and $\psi$-types over finite sets respectively.  Then notice that 
 \[
  (\gamma_\varphi \wedge \gamma_\psi)(\overline{y} \concat \overline{z}; \overline{w}_0 \concat \overline{v}_0, ..., \overline{w}_{n-1} \concat \overline{v}_{n-1}) = \gamma_\varphi(\overline{y}; \overline{w}_0, ..., \overline{w}_{n-1}) \wedge \gamma_\psi(\overline{z}; \overline{v}_0, ..., \overline{v}_{n-1})
 \]
 is a uniform definition of $(\varphi \wedge \psi)$-types over finite sets and notice that $\neg \gamma_\varphi(\overline{y}; \overline{w}_0, ..., \overline{w}_{n-1})$ is a uniform definition of $(\neg \varphi)$-types over finite sets.  This yields the desired conclusion by induction on formula complexity.
\end{proof}

The next proposition follows by definition and by the characterization of stability in terms of uniform definability of types (see Theorem II.2.12 in \cite{shelah}).

\begin{prop}\label{Prop_StableUDTFS}
 If $\varphi(\overline{x}; \overline{y})$ is a stable formula, then $\varphi$ has UDTFS.  Thus, stable theories have UDTFS.
\end{prop}

The next proposition puts the property of UDTFS for formulas between stability and dependence.

\begin{prop}\label{Prop_UDTFSDep}
 If $\varphi(\overline{x}; \overline{y})$ has UDTFS, then $\varphi$ is dependent.
\end{prop}

\begin{proof}
 If $\varphi(\overline{x}; \overline{y})$ has UDTFS, then it certainly has uniform definability of types over finite indiscernible sequences.  Therefore, it satisfies the right-hand side of Theorem \ref{Thm_IndiscernibleSequenceChar} (ii), hence is dependent.
\end{proof}

It is still open whether or not dependent formulas have UDTFS or even if dependent theories have UDTFS (see Open Question \ref{Conj_DepUDTFS} below).

The next lemma shows that we do not actually need a single uniform definition of $\varphi$-types over finite sets, but it suffices to have a fixed finite number of uniform definitions of $\varphi$-types over finite sets.  This simplifies arguments showing that formulas and theories have UDTFS.  This is essentially due to Shelah in the proof of Theorem II.2.12 (1) in \cite{shelah}, where he shows it for the standard definability of types in the same manner.

\begin{lem}\label{Lemma_ManyOneUDTFS}
 Fix $\varphi(\overline{x}; \overline{y})$ a partitioned formula.  If there exists $\{ \psi_\ell(\overline{y}; \overline{z}_0, ... \overline{z}_{N_\ell - 1}) : \ell < L \}$ a finite collection of formulas such that, for all finite $B$ with $|B| \ge 2$ and for all $p(\overline{x}) \in S_{\varphi}(B)$, there exists $\ell < L$ and $\overline{c}_0, ..., \overline{c}_{N_\ell-1} \in B$ such that $\psi_\ell(\overline{y}; \overline{c}_0, ..., \overline{c}_{N_\ell-1})$ defines $p(\overline{x})$, then $\varphi$ has UDTFS.
\end{lem}

\begin{proof}
 For simplicity, assume $N_\ell = N$ for all $\ell < L$.  Consider the following formula:
 \[
  \psi(\overline{y}; \overline{z}_0, ..., \overline{z}_{N-1}, \overline{w}, \overline{v}_0, ..., \overline{v}_{L-1}) = \bigwedge_{\ell < L} (\overline{w} = \overline{v}_\ell \rightarrow \psi_\ell(\overline{y}; \overline{z}_0, ..., \overline{z}_{N-1})).
 \]
 Then, fix any finite set of $\leng(\overline{y})$-tuples $B$ with $|B| \ge 2$ and fix any $p(\overline{x}) \in S_{\varphi}(B)$.  By hypothesis, there exists $\ell < L$ and $\overline{c}_0, ..., \overline{c}_{N-1} \in B$ so that $\psi_\ell(\overline{y}; \overline{c}_0, ..., \overline{c}_{N-1})$ defines $p$.  Fix any $\overline{b} \neq \overline{b}'$ from $B$ (here we use the hypothesis that $|B| \ge 2$) and let $\overline{b}_i = \overline{b}$ unless $i = \ell$, in which case let $\overline{b}_\ell = \overline{b}'$.  Then the following formula defines $p$:
 \[
  \psi(\overline{y}; \overline{c}_0, ..., \overline{c}_{N-1}, \overline{b}', \overline{b}_0, ..., \overline{b}_{L-1}).
 \]
 Therefore, we see that $\psi$ is a uniform definition of $\varphi$-types over finite sets, so $\varphi$ has UDTFS.
\end{proof}

We exhibit another lemma that reduces the difficulty in showing that a theory has UDTFS.

\begin{lem}[Sufficiency of a single variable]\label{Lemma_SuffofSingleVar}
 If $T$ is a theory such that all formulas $\varphi(x; \overline{y})$ have UDTFS (where $x$ is a singleton), then $T$ has UDTFS.
\end{lem}

\begin{proof}
 This proceeds by induction on $\leng(\overline{x})$, with the base case taken care of by assumption.

 Consider $\leng(\overline{x}) = n$ and assume that all formulas with less than $n$ free variables have UDTFS.  Then consider $\varphi(\overline{x}; \overline{y})$.  We construct a uniform definition of $\varphi$-types over finite sets.

 Consider the repartitioned $\hat{\varphi}(x_1 \concat x_2 \concat ... \concat x_{n-1}; x_n \concat \overline{y}) = \varphi(\overline{x}; \overline{y})$.  Now this has only $n-1$ free variables, so by induction, there exists a uniform definition of $\hat{\varphi}$-types over finite sets.  Let that formula be $\psi(x_n \concat \overline{y}; w_0 \concat \overline{z}_0, ..., w_{k-1} \concat \overline{z}_{k-1})$ (where $\leng(\overline{z}_i) = \leng(\overline{y})$).  Now let $\psi'(x_n \concat \overline{y}; \overline{z}_0, ..., \overline{z}_{k-1}) = \psi(x_n \concat \overline{y}; x_n \concat \overline{z}_0, ..., x_n \concat \overline{z}_{k-1})$, where we replace each $w_i$ with $x_n$.  Consider the repartitioned $\hat{\psi'}(x_n; \overline{y} \concat \overline{z}_0 \concat ... \concat \overline{z}_{k-1}) = \psi'(x_n \concat \overline{y}; \overline{z}_0, ..., \overline{z}_{k-1})$.  Again, this has only one free variable, so by hypothesis, there exists a uniform definition of $\hat{\psi'}$-types over finite sets, say $\gamma(\overline{y} \concat \overline{z}_0 \concat ... \concat \overline{z}_{k-1}; \overline{\mathbf{w}})$ (where $\overline{\mathbf{w}}$ has length a multiple of $(k+1) \cdot \leng(\overline{y})$, as the dependent variables $\overline{y} \concat \overline{z}_0 \concat ... \concat \overline{z}_{k-1}$ have length $(k+1) \cdot \leng(\overline{y})$).  Consider the repartitioned $\hat{\gamma}(\overline{y}; \overline{z}_0, ..., \overline{z}_{k-1}, \overline{\mathbf{w}}) = \gamma(\overline{y} \concat \overline{z}_0 \concat ... \concat \overline{z}_{k-1}; \overline{\mathbf{w}})$.  I claim that this is a uniform definition of $\varphi$-types over finite sets, as desired.

 Take $B$ a finite set of $\leng(\overline{y})$-tuples from $\mathfrak{C}$ and $\overline{a} = (a_1, ..., a_n) \in \mathfrak{C}^n$, and consider the $\varphi$-type $p(\overline{x}) = \tp_{\varphi}(\overline{a}/B)$.  Now consider the $\hat{\varphi}$-type $\hat{p}(x_1, ... , x_{n-1}) = \tp_{\hat{\varphi}}(a_1, ..., a_{n-1} / a_n \concat B)$ (where we define $a_n \concat B = \{ a_n \concat \overline{b} : \overline{b} \in B \}$).  As $\psi$ is a finite definition of $\hat{\varphi}$-types, there exists some $\overline{c}_0, ..., \overline{c}_{k-1} \in a_n \concat B$ such that:

 \begin{center}
  For all $\overline{b} \in B$, $\hat{\varphi}(x_1, ..., x_{n-1}; a_n \concat \overline{b}) \in \hat{p}$ if and only if $\models \psi(a_n \concat \overline{b}; \overline{c}_0, ..., \overline{c}_{k-1})$.
 \end{center}

 But we notice that $\overline{c}_i \in a_n \concat B$ means that $\overline{c}_i = a_n \concat \overline{d}_i$ for some $\overline{d}_i \in B$.  With this substitution, we see that for all $\overline{b} \in B$, $\hat{\varphi}(x_1, ..., x_{n-1}; a_n \concat \overline{b}) \in \hat{p}$ if and only if $\models \psi'(a_n \concat \overline{b}; \overline{d}_0, ..., \overline{d}_{k-1})$.  Now consider $q(x_n) = \tp_{\hat{\psi'}}(a_n / B^{k+1})$ (note that the use of $\psi'$ was to remove the need to consider types over $a_n \concat B$, and instead consider them over $B^{k+1}$).  Using the uniform definition of $\hat{\psi'}$-types over finite sets, $\gamma$, we get that there exists a tuple $\overline{\mathbf{e}} \in B^{\leng(\overline{\mathbf{w}}) / \leng(\overline{y})}$ such that:

 \begin{center}
  For all $\overline{b} \in B$, $\models \psi'(a_n; \overline{b} \concat \overline{d}_0 \concat ... \concat \overline{d}_{k-1})$ if and only if $\models \gamma(\overline{b} \concat \overline{d}_0 \concat ... \concat \overline{d}_{k-1}; \overline{\mathbf{e}})$.
 \end{center}

 But then we have that, for all $\overline{b} \in B$:

 $\varphi(x_1, ..., x_{n-1}, x_n; \overline{b}) \in p$ iff. $\hat{\varphi}(x_1, ..., x_{n-1}; a_n \concat \overline{b}) \in \hat{p}$ iff. $\models \psi(a_n \concat \overline{b}; \overline{c}_0, ..., \overline{c}_{k-1})$ iff. $\models \psi'(a_n \concat \overline{b}; \overline{d}_0, ..., \overline{d}_{k-1})$ iff. $\models \gamma(\overline{b} \concat \overline{d}_0 \concat ... \concat \overline{d}_{k-1}; \overline{\mathbf{e}})$ iff. $\models \hat{\gamma}(\overline{b}; \overline{d}_0, ..., \overline{d}_{k-1}, \overline{\mathbf{e}})$.

 That is, $\hat{\gamma}$ is a uniform definition of $\varphi$-types over finite sets.
\end{proof}

The sufficiency of a single variable lemma is used to prove UDTFS for a theory whose unary formulas are well understood.  For example, one can use it to show that weakly o-minimal theories have UDTFS.  A linearly ordered structure $(M; <, ...)$ is \textit{weakly o-minimal} if all parameter-definable subsets of $M$ are a union of finitely many $<$-convex sets.  A theory $T$ is \textit{weakly o-minimal} if all models of $T$ are weakly o-minimal.  This next proposition was known to Johnson and Laskowski.  It follows as a corollary of Theorem \ref{Thm_dpMinUDTFS} below.

\begin{prop}\label{Prop_wominUDTFS}
 If $T$ is weakly o-minimal, then $T$ has UDTFS.
\end{prop}

The goal of the next section is to generalize Proposition \ref{Prop_wominUDTFS} by showing that all dp-minimal theories have UDTFS.


\section{dp-Minimal Theories}

Fix $T$ a complete theory and let $\mathfrak{C}$ be a monster model of $T$.  We aim to prove the following theorem:

\begin{thm}\label{Thm_dpMinUDTFS}
 If $T$ is dp-minimal, then $T$ has UDTFS.
\end{thm}

We begin by defining ICT-patterns and dp-minimality.  This definition is taken from \cite{DGL}, for example.

\begin{defn}\label{Def_ICTpattern}
 An \textit{ICT pattern} is a pair of formulas $\varphi(x; \overline{y})$ and $\psi(x; \overline{z})$ with a single free variable $x$ and sequences $\langle \overline{b}_i : i < \omega \rangle$ and $\langle \overline{c}_j : j < \omega \rangle$ such that, for all $\ell, k < \omega$, the following partial type is consistent:
 \[
  \{ \varphi(x; \overline{b}_\ell), \psi(x; \overline{c}_k) \} \cup \{ \neg \varphi(x; \overline{b}_i) : i < \omega, i \neq \ell \} \cup \{ \neg \psi(x; \overline{c}_j) : j < \omega, j \neq k \}.
 \]
 A theory $T$ is \textit{dp-minimal} if there exists no ICT pattern.
\end{defn}

Consider instead TP-patterns, defined as follows:

\begin{defn}
 A \textit{TP-pattern} is a formula $\varphi(x; \overline{y})$ with a single free variable $x$ and a sequence of $\leng(\overline{y})$-tuples $\langle \overline{b}_i : i < \omega \rangle$ such that, for all $\ell < k < \omega$, the following formula holds:
 \[
  \exists x \left( \varphi(x; \overline{b}_k) \wedge \varphi(x; \overline{b}_{\ell}) \wedge \bigwedge_{i < k, i \neq \ell} \neg \varphi(x; \overline{b}_i) \right).
 \]
\end{defn}

We show that having a TP-pattern is equivalent to having an ICT pattern.

\begin{prop}\label{Prop_dpMintpMin}
 If $T$ is dp-minimal if and only if $T$ has no TP-pattern
\end{prop}

\begin{proof}
 Suppose that $T$ has a TP-pattern and fix $\varphi(x; \overline{y})$ and $\langle \overline{b}_i : i < \omega \rangle$ such a TP-pattern.  Let $\psi(x; \overline{y}_0, \overline{y}_1) = \varphi(x; \overline{y}_0) \leftrightarrow \varphi(x; \overline{y}_1)$ and let $K$ be any positive integer.  By Ramsey's theorem, we may assume that $\langle \overline{b}_i : i < \omega \rangle$ is $\Delta$-indiscernible, where $\Delta = \Delta_{4K}$ (as defined in \eqref{Eq_defofDeltaN}).  By definition of a TP-pattern, we know that the following is consistent:
 \[
  \{ \neg \varphi(x; \overline{b}_i) : i < 2K \} \cup \{ \varphi(x; \overline{b}_{2K}) \} \cup \{ \neg \varphi(x; \overline{b}_i) : 2K < i \le 6K \} \cup \{ \varphi(x; \overline{b}_{6K+1}) \}.
 \]
 Let $a$ realize this type.  By trimming down the sequence, we may assume that the truth value of $\varphi(a; \overline{b}_i)$ is constant for all $i > 6K+1$.  We therefore have that the following is consistent, witnessed by $a$:
 \begin{align*}
  \{ \psi(x; \overline{b}_{2i}, \overline{b}_{2i+1}) : i < K \} \cup \{ \neg \psi(x; \overline{b}_{2K}, \overline{b}_{2K+1}) \} \cup \{ \psi(x; \overline{b}_{2i}, \overline{b}_{2i+1}) : K < i < 3K \} \\ \cup \{ \neg \psi(x; \overline{b}_{6K}, \overline{b}_{6K+1}) \} \cup \{ \psi(x; \overline{b}_{2i}, \overline{b}_{2i+1}) : 3K < i < 4K \}.
 \end{align*}
 By $\Delta$-indiscernibility, we therefore have that, for all $\ell < K$ and $K \le k < 2K$, the following is consistent:
 \begin{align*}
  \{ \psi(x; \overline{b}_{2i}, \overline{b}_{2i+1}) : i < K, i \neq \ell \} \cup \{ \neg \psi(x; \overline{b}_{2\ell}, \overline{b}_{2\ell+1}) \} \cup \\ \{ \psi(x; \overline{b}_{2i}, \overline{b}_{2i+1}) : K < i < 2K, i \neq k \} \cup \{ \neg \psi(x; \overline{b}_{2k}, \overline{b}_{2k+1}) \}.
 \end{align*}
 Since $K$ was arbitrary, by compactness there exists $\overline{c}_i, \overline{d}_j \in \mathfrak{C}^{\leng(\overline{y})}$ for all $i,j < \omega$ such that, for all $\ell, k < \omega$,
 \begin{align*}
  \{ \psi(x; \overline{c}_{2i}, \overline{c}_{2i+1}) : i < \omega, i \neq \ell \} \cup \{ \neg \psi(x; \overline{c}_{2\ell}, \overline{c}_{2\ell+1}) \} \cup \\ \{ \psi(x; \overline{d}_{2j}, \overline{d}_{2j+1}) : j < \omega, j \neq k \} \cup \{ \neg \psi(x; \overline{d}_{2k}, \overline{d}_{2k+1}) \}
 \end{align*}
 is consistent.  Then $\neg \psi$, $\neg \psi$, $\langle \overline{c}_{2i} \concat \overline{c}_{2i+1} : i < \omega \rangle$, and $\langle \overline{d}_{2j} \concat \overline{d}_{2j+1} : j < \omega \rangle$ form an ICT-pattern.  Thus, $T$ is not dp-minimal.

 Conversely, suppose that $T$ is not dp-minimal and let $\varphi(x; \overline{y})$, $\psi(x; \overline{z})$, $\langle \overline{b}_i : i < \omega \rangle$, and $\langle \overline{c}_j : j < \omega \rangle$ be an ICT-pattern.  Then let $\theta(x; \overline{y}, \overline{z}) = \neg (\varphi(x; \overline{y}) \leftrightarrow \psi(x; \overline{z}))$.  It is clear that $\theta$ along with $\langle \overline{b}_i \concat \overline{c}_i : i < \omega \rangle$ form a TP-pattern.
\end{proof}

Assume $T$ is dp-minimal, hence has no TP-pattern, and fix $\varphi(x; \overline{y})$ any partitioned formula.

\begin{rem}\label{Rem_tpMinBound}
 By compactness, there exists $K < \omega$ such that, for all sequences $\langle \overline{b}_i : i < K \rangle$ with $\overline{b}_i \in \mathfrak{C}^{\leng(\overline{y})}$ for all $i$, we have that

 \begin{itemize}
  \item [(i)] There exists $\ell < k < K$ such that $\models \neg \exists x \left(\varphi(x; \overline{b}_k) \wedge \varphi(x; \overline{b}_{\ell}) \wedge \bigwedge_{i < k, i \neq \ell} \neg \varphi(x; \overline{b}_i) \right)$, and
  \item [(ii)] There exists $\ell < k < K$ such that $\models \neg \exists x \left(\neg \varphi(x; \overline{b}_k) \wedge \neg \varphi(x; \overline{b}_{\ell}) \wedge \bigwedge_{i < k, i \neq \ell} \varphi(x; \overline{b}_i) \right)$.
 \end{itemize}
\end{rem}

Fix $B \subseteq \mathfrak{C}^{\leng(\overline{y})}$ finite and $p \in S_\varphi(B)$ any type.  Let $\delta : B \rightarrow 2$ be the function associated to $p$ (i.e. $\varphi(x; \overline{b})^{\delta(\overline{b})} \in p(x)$ for all $\overline{b} \in B$).  We first begin by defining what it means for a small $\varphi$-type to decide a formula.

\begin{defn}\label{Def_Decides}
 For any $\varphi$-type $q(x)$ and any $\overline{b} \in B$, we say that \textit{$q$ decides $\varphi(x; \overline{b})$} if either $q(x) \vdash \varphi(x; \overline{b})$ or $q(x) \vdash \neg \varphi(x; \overline{b})$.  We say that \textit{$q$ decides $\varphi(x; \overline{b})$ correctly} if $q(x) \vdash \varphi(x; \overline{b})^{\delta(\overline{b})}$ (i.e. $q$ decides $\varphi(x; \overline{b})$ and it implies the instance of $\pm \varphi(x; \overline{b})$ that is actually in $p(x)$).
\end{defn}

Notice that $q$ need not be a $\varphi$-subtype of $p$.  When it is, the following lemma is immediate since $p$ is consistent:

\begin{lem}\label{Lem_ObviousDecide}
 For any $\varphi$-subtype $q(x) \subseteq p(x)$ and any $\overline{b} \in B$, if $q$ decides $\varphi(x; \overline{b})$ then it does so correctly.
\end{lem}

For any subsets $B_1 \subseteq B_0 \subseteq B$, recall the definition of $p_{B_0}$ and $p_{B_0, B_1}$ from the introduction ($p_{B_0}$ is the restriction of $p$ to $B_0$ and $p_{B_0, B_1}$ is $p_{B_0}$ perturbed by $B_1$).

\begin{defn}\label{Def_nDecides}
 Fix $B_0 \subseteq B$ and $\overline{b} \in B$.  We say that \textit{$B_0$ $*$-decides $\varphi(x; \overline{b})$} if $p_{B_0}$ decides $\varphi(x; \overline{b})$ or there exists $\overline{b}' \in B_0$ such that $p_{B_0, \{ \overline{b}' \}}$ is consistent and decides $\varphi(x; \overline{b})$.  We say that \textit{$B_0$ $*$-decides $\varphi(x; \overline{b})$ correctly} if $B_0$ $*$-decides $\varphi(x; \overline{b})$ and we have that either
 \begin{itemize}
  \item [(i)] $p_{B_0}$ decides $\varphi(x; \overline{b})$, or
  \item [(ii)] for all $\overline{b}' \in B_0$ such that $p_{B_0, \{ \overline{b}' \}}$ is consistent and decides $\varphi(x; \overline{b})$, $p_{B_0, \{ \overline{b}' \}}$ decides $\varphi(x; \overline{b})$ correctly.
 \end{itemize}
\end{defn}

So $B_0$ $*$-decides $\varphi(x; \overline{b})$ if there is a perturbation of $p_{B_0}$ of size at most one that decides $\varphi(x; \overline{b})$.  By Lemma \ref{Lem_ObviousDecide}, if $p_{B_0}$ decides $\varphi(x; \overline{b})$ then it does so correctly.  Therefore, the only way $B_0$ would $*$-decide $\varphi(x; \overline{b})$ incorrectly is if $p_{B_0}$ does not decide $\varphi(x; \overline{b})$ and, for some $\overline{b}' \in B_0$, $p_{ B_0, \{ \overline{b}' \} }$ is consistent and decides $\varphi(x; \overline{b})$ incorrectly (i.e. $p_{ B_0, \{ \overline{b}' \} } (x) \vdash \neg \varphi(x; \overline{b})^{\delta(\overline{b})}$).

The concept of $*$-decides captures one possible way of constructing an algorithm to define the $\varphi$-type, $p$.  If we can construct, in a uniform manner, a small collection of small subsets of $B$ which, when chosen in a certain order, $*$-decides $\varphi(x; \overline{b})$ correctly for all $\overline{b} \in B$, we can get a uniform definition of $\varphi$-types over finite sets.  We now show how to construct such collection.

To aid notation and cast these ideas in terms of an ordering, we define a quasi-ordering on the powerset of $B$, $\mathcal{P}(B)$, as follows:

\begin{center}
 For $B_0, B_1 \in \mathcal{P}(B)$, let $B_0 \le_p B_1$ if $p_{B_0} (x) \vdash p_{B_1} (x)$.
\end{center}

We say that $B_0$ is \textit{$p$-equivalent} to $B_1$, denoted $B_0 \equiv_p B_1$, if $B_0 \le_p B_1$ and $B_1 \le_p B_0$.  Notice that $\equiv_p$ is an equivalence relation of $\mathcal{P}(B)$ and $\le_p$ is a partial ordering on $\mathcal{P}(B) / \equiv_p$.  Say that $B_0 <_p B_1$ if $B_0 \le_p B_1$ but $B_1 \not\le_p B_0$ (i.e. $B_0 \not\equiv_p B_1$).  For completeness, set $p_{\emptyset}(x) = (x = x)$.  So, $B_0 \le_p \emptyset$ for all $B_0 \in \mathcal{B}(B)$ and $\emptyset \le_p B_0$ if and only if $p_{B_0}$ is realized by all elements of $\mathfrak{C}$.  The following lemma is immediate from the definitions:

\begin{lem}\label{Lem_QuasiOrderingp}
 For the quasi-ordering $\le_p$, the following hold:

 \begin{itemize}
  \item [(i)] For all $B_0, B_1 \in \mathcal{P}(B)$, $B_0 \subseteq B_1$ implies that $B_1 \le_p B_0$.
  \item [(ii)] For all $B_0, B_1 \in \mathcal{P}(B)$ and all $B'_0 \subseteq B_0$, $B_0 \le_p B_1$ if and only if $B_0 \le_p B_1 \cup B'_0$.
  \item [(iii)] If $\mathcal{B} \subseteq \mathcal{P}(B)$ and $B_1 \in \mathcal{B}$, then there exists a $B_0 \le_p B_1$ such that $B_0 \in \mathcal{B}$ and, for all other $B_2 \in \mathcal{B}$, $B_2 \le_p B_0$ implies that $B_2 \equiv_p B_0$ (we call such $B_0$ $\le_p$-minimal elements of $\mathcal{B}$).
 \end{itemize}
\end{lem}

Notice that (iii) holds because $B$, hence $\mathcal{B}$, is finite.  This is our main use of finiteness and the only obstacle for showing general uniform definability for dp-minimal theories.  A great deal of milage can be obtained by using $\le_p$-minimal elements.  Let $\mathcal{B}$ be any non-empty set of subsets of $B$.  Consider the following lemma about correct decisions that uses $\le_p$-minimality:

\begin{lem}\label{Lem_MakingCorrectDecisions}
 Fix $\overline{b} \in B$, $B_0 \in \mathcal{B}$ $\le_p$-minimal in $\mathcal{B}$, and $\overline{b}' \in B_0$.  If

 \begin{itemize}
  \item [(i)] $p_{B_0, \{ \overline{b}' \}}$ is consistent and decides $\varphi(x; \overline{b})$,
  \item [(ii)] $p_{B_0}$ does not decide $\varphi(x; \overline{b})$, and
  \item [(iii)] there exists $B_1 \le_p B_0 - \{ \overline{b}' \}$ such that $B_1 \cup \{ \overline{b} \} \in \mathcal{B}$
 \end{itemize}

 then $p_{B_0, \{ \overline{b}' \}}$ decides $\varphi(x; \overline{b})$ correctly.
\end{lem}

\begin{proof}
 Since $p_{B_0}$ does not decide $\varphi(x; \overline{b})$, by definition of $\le_p$, we have that $B_0 \not\le_p \{ \overline{b} \}$, hence, by Lemma \ref{Lem_QuasiOrderingp} (ii),
 \begin{equation}\label{Eq_NotEquivEq}
  B_0 \not\equiv_p \left( B_0 - \{ \overline{b}' \} \right) \cup \{ \overline{b} \}.
 \end{equation}
 Now, by means of contradiction, suppose that $p_{B_0, \{ \overline{b}' \}}$ decides $\varphi(x; \overline{b})$ incorrectly.  That is, suppose that
 \[
  p_{(B_0 - \{ \overline{b}' \})}(x) \cup \left\{ \neg \varphi(x; \overline{b}')^{\delta(\overline{b}')} \right\} \vdash \neg \varphi(x; \overline{b})^{\delta(\overline{b})}.
 \]
 By contrapositive, we get that $p_{(B_0 - \{ \overline{b}' \}) \cup \{ \overline{b} \}} \vdash p_{\{ \overline{b}' \}}$, hence $( B_0 - \{ \overline{b}' \} ) \cup \{ \overline{b} \} \le_p \{ \overline{b}' \}$.  By Lemma \ref{Lem_QuasiOrderingp} (ii), we get that $( B_0 - \{ \overline{b}' \} ) \cup \{ \overline{b} \} \le_p B_0$.  However, by \eqref{Eq_NotEquivEq}, we get that this is strict, so
 \[
  \left( B_0 - \{ \overline{b}' \} \right) \cup \{ \overline{b} \} <_p B_0.
 \]
 Now, using (iii) we note that 
 \[
  B_1 \cup \{ \overline{b} \} \le_p \left( B_0 - \{ \overline{b}' \} \right) \cup \{ \overline{b} \} <_p B_0
 \]
 and $B_1 \cup \{ \overline{b} \} \in \mathcal{B}$.  This contradicts the $\le_p$-minimality of $B_0$ in $\mathcal{B}$.
\end{proof}

So to get correct decisions, we need only find a $\mathcal{B}$ that has good closure properties.  We now construct such a $\mathcal{B}$.  First, notice that everything worked out above for subsets of $B$ translates to sequences in $B$ by considering the images of those sequences.  With this in mind, we define $\mathcal{B}_n$, a set of sequences from $B$ of length $n$, inductively as follows:

For $n=1$, let $\mathcal{B}_1 = \{ \langle \overline{b} \rangle : \overline{b} \in B \text{ and } \varphi(x; \overline{b})^t \text{ is consistent for both } t < 2 \}$ (i.e. for all $\overline{b} \in B$, $\langle \overline{b} \rangle \in \mathcal{B}_1$ if and only if $\emptyset \not\equiv_p \langle \overline{b} \rangle$ if and only if $p_{ \langle \overline{b} \rangle, \langle \overline{b} \rangle }$ is consistent).  For $n > 1$, let
\[
 \mathcal{B}_n = \{ \beta \concat \langle \overline{b} \rangle : \overline{b} \in B, \beta \in \mathcal{B}_{n-1}, \text{ and } \beta \text{ does not $*$-decide } \varphi(x; \overline{b}) \}.
\]
We start off with a basic Lemma about elements of $\mathcal{B}_n$.

\begin{lem}\label{Lem_BasicBn}
 If $\beta = \langle \overline{b}_0, ..., \overline{b}_{n-1} \rangle \in \mathcal{B}_n$, then the following hold:

 \begin{itemize}
  \item [(i)] For all $\ell < n$, $p_{\beta, \langle \overline{b}_\ell \rangle}$ is consistent.
  \item [(ii)] For all $\ell < k < n$, $p_{\langle \overline{b}_i : i \le k \rangle, \langle \overline{b}_\ell, \overline{b}_k \rangle}$ is consistent.
  \item [(iii)] For all $k \le n$ and all subsequences $\beta_0 \subseteq \beta$ of length $k$, $\beta_0 \in \mathcal{B}_k$.
  \item [(iv)] $\left| S_\varphi( \{ \overline{b}_i : i < n \} ) \right| > \frac{n(n+1)}{2}$.
  \item [(v)] For all $\overline{b} \in B$, if $\beta$ does not $*$-decide $\varphi(x; \overline{b})$, then $\beta \concat \langle \overline{b} \rangle \in \mathcal{B}_{n+1}$.
 \end{itemize}
\end{lem}

\begin{proof}
 (i) and (ii): This follows by induction on $n$ and is clear for $n = 1$.  Fix $n > 1$ and suppose that both (i) and (ii) hold below $n$.  Let $\beta' = \langle \overline{b}_i : i < n - 1 \rangle$, the initial segment of $\beta$ of length $n - 1$.  Thus we have that $p_{\beta', \langle \overline{b}_\ell \rangle}$ is consistent for all $\ell < n - 1$ and $p_{\langle \overline{b}_i : i \le k \rangle, \langle \overline{b}_\ell, \overline{b}_k \rangle}$ is consistent for all $\ell < k < n - 1$.  However, $\beta'$ does not $*$-decide $\varphi(x; \overline{b}_{n-1})$ by definition of $\mathcal{B}_n$.  Therefore, it is clear that $p_{\beta, \langle \overline{b}_\ell \rangle}$ is consistent and $p_{\beta, \langle \overline{b}_\ell, \overline{b}_{n-1} \rangle}$ is consistent for all $\ell < n$.  Therefore (i) and (ii) hold on $n$.

 (iii): Fix $k < n$ and any subsequence $\beta_0 \subseteq \beta$ and let $\beta_0 = \langle \overline{b}_{i(0)}, ..., \overline{b}_{i(k-1)} \rangle$.  We now show $\beta_0 \in \mathcal{B}_k$ by induction on $k$.  If $k=1$, then (i) implies that $p_{ \langle \overline{b}_{i(0)} \rangle, \langle \overline{b}_{i(0)} \rangle }$ is consistent, hence $\beta_0 \in \mathcal{B}_1$.  Now, assuming by induction that $\langle \overline{b}_{i(0)}, ..., \overline{b}_{i(k-2)} \rangle \in \mathcal{B}_{k-1}$, it suffices to show that $\langle \overline{b}_{i(0)}, ..., \overline{b}_{i(k-2)} \rangle$ does not $*$-decide $\overline{b}_{i(k-1)}$.  This follows from conditions (i) and (ii) (``does not $*$-decide'' is implied by some of the types in (i) and (ii) being consistent).

 (iv): This follows from (i) and (ii) by counting (note that $p_\beta$ is also consistent, giving the strict inequality).

 (v): This follows by definition of $\mathcal{B}_n$.
\end{proof}

We note that the non-existence of a TP-pattern gives us a bound on the $n < \omega$ such that $\mathcal{B}_n \neq \emptyset$.

\begin{lem}\label{Lem_tpMinBoundBn}
 Given $K < \omega$ as in Remark \ref{Rem_tpMinBound}, $\mathcal{B}_{2K-1} = \emptyset$.
\end{lem}

\begin{proof}
 By means of contradiction, fix $\beta \in \mathcal{B}_{2K-1}$.  By pigeon-hole principle, there exists $t < 2$ and a subsequence $\langle \overline{b}_0, ..., \overline{b}_{K-1} \rangle = \beta_0 \subseteq \beta$ of length $K$ such that, for all $i < K$, $\delta(\overline{b}_i) = t$ (i.e. $\delta$, hence $p$, is constant on $\beta_0$).  By Lemma \ref{Lem_BasicBn} (iii), $\beta_0 \in \mathcal{B}_K$.  Therefore, by Lemma \ref{Lem_BasicBn} (ii), for all $\ell < k < K$,
 \[
  p_{ \{ \overline{b}_i : i \le k \}, \{ \overline{b}_\ell, \overline{b}_k \} } (x)
 \]
 is consistent.  However, since $p$ is constant on $\beta_0$, we have that, for all $\ell < k < K$,
 \[
  \{ \varphi(x, \overline{b}_i)^t : i < k, i \neq \ell \} \cup \{ \neg \varphi(x; \overline{b}_\ell)^t, \neg \varphi(x; \overline{b}_k)^t \}
 \]
 is consistent.  This contradicts our choice of $K$ as in Remark \ref{Rem_tpMinBound} (i.e. it contradicts the non-existence of a TP-pattern of length $K$).
\end{proof}

We now define, for each $n$ and each $\beta$ a $\le_p$-minimal element of $\mathcal{B}_n$, a set $H(\beta)$ of non-empty sequences inductively as follows:

For $n = 1$, fix $\beta \in \mathcal{B}_1$ that is $\le_p$-minimal in $\mathcal{B}_1$ and let $H(\beta) = \{ \beta \}$.  For $n > 1$, fix any $\beta \in \mathcal{B}_n$ that is $\le_p$-minimal in $\mathcal{B}_n$.  Let $\beta = \langle \overline{b}_0, ..., \overline{b}_{n-1} \rangle$ and, for each $i < n$, let $\beta_i$ be the subsequence of $\beta$ given by
\[
 \beta_i = \langle \overline{b}_0, ..., \overline{b}_{i-1}, \overline{b}_{i+1}, ..., \overline{b}_{n-1} \rangle
\]
 (i.e. by removing the $i$th element from $\beta$).  By Lemma \ref{Lem_BasicBn} (iii), $\beta_i \in \mathcal{B}_{n-1}$.  Therefore, by Lemma \ref{Lem_QuasiOrderingp} (iii), there exists a $\le_p$-minimal element of $\mathcal{B}_{n-1}$, $\beta'_i$, such that $\beta'_i \le_p \beta_i$.  Fix any choice of $\beta'_i$ for each $i < n$.  Finally, let
\[
 H(\beta) = \bigcup_{i < n} H(\beta'_i) \cup \{ \beta \}.
\]
This defines $H$ on all $\le_p$-minimal elements of $\mathcal{B}_n$ for each $n$, as desired.  Note that $H(\beta)$ has size a function of $\leng(\beta)$ (bounded, for example, by $(\leng(\beta) + 1)!$) and each element of $H(\beta)$ is a non-empty sequence of length at most $\leng(\beta)$.  We now show that elements of $H(\beta)$, when chosen in a particular manner, can correctly $*$-decide $\varphi(x; \overline{b})$.

\begin{lem}\label{Lem_CorrectDecisionsonH}
 Fix $\overline{b} \in B$, $n < \omega$, and $\beta \in \mathcal{B}_n$ $\le_p$-minimal.  Let $k \le n$ be minimal such that there exists $\beta' \in H(\beta)$ with $\leng(\beta') = k$ and $\beta'$ $*$-decides $\varphi(x; \overline{b})$.  Then any such $\beta'$ $*$-decides $\varphi(x; \overline{b})$ correctly.
\end{lem}

\begin{proof}
 By induction on $n$.  If $n = 1$, then we need only check that if $\beta \in \mathcal{B}_1$ is $\le_p$-minimal and $\beta$ $*$-decides $\varphi(x; \overline{b})$, then it does so correctly.  If $p_\beta$ decides $\varphi(x; \overline{b})$, then it does so correctly by Lemma \ref{Lem_ObviousDecide}.  So it remains to show that if $p_\beta$ does not decide $\varphi(x; \overline{b})$ and $p_{\beta, \beta}$ decides $\varphi(x; \overline{b})$, then $p_{\beta, \beta}$ decides $\varphi(x; \overline{b})$ correctly.  This follows from Lemma \ref{Lem_MakingCorrectDecisions} (note that, since $p_\beta$ does not decide $\varphi(x; \overline{b})$, $\langle \overline{b} \rangle \in \mathcal{B}_1$).

 Suppose now that $n > 1$.  Let $k \le n$ be minimal such that there exists $\beta' \in H(\beta)$ with $\leng(\beta') = k$ and $\beta'$ $*$-decides $\varphi(x; \overline{b})$.  By induction, we may assume that $k = n$ and $\beta' = \beta$.  By Lemma \ref{Lem_ObviousDecide}, we may assume that $p_\beta$ does not decide $\varphi(x; \overline{b})$.  So assume that $p_{\beta, \langle \overline{b}_\ell \rangle}$ decides $\varphi(x; \overline{b})$, where we let $\beta = \langle \overline{b}_0, ..., \overline{b}_{n-1} \rangle$.  Consider $\beta'_\ell \le_p \beta_\ell$ as defined above.  We have that $\beta'_\ell \in H(\beta)$ and $\leng(\beta'_\ell) = n - 1 < n$.  By minimality, $\beta'_\ell$ does not $*$-decide $\varphi(x; \overline{b})$.  By Lemma \ref{Lem_BasicBn} (v), we know that $\beta'_\ell \concat \langle \overline{b} \rangle \in \mathcal{B}_n$.  These are exactly the conditions needed in Lemma \ref{Lem_MakingCorrectDecisions}.  Therefore, we get that $p_{\beta, \langle \overline{b}_\ell \rangle}$ decides $\varphi(x; \overline{b})$ correctly.  Since $\ell < n$ was arbitrary such that $p_{\beta, \langle \overline{b}_\ell \rangle}$ decides $\varphi(x; \overline{b})$, we see that $\beta$ $*$-decides $\varphi(x; \overline{b})$ correctly.
\end{proof}

This correct decision process is enough to prove UDTFS.

\begin{proof}[Proof of Theorem \ref{Thm_dpMinUDTFS}]
 By sufficiency of a single variable, it suffices to show that formulas of the form $\varphi(x; \overline{y})$ have UDTFS.  Fix such a formula and let $K < \omega$ be given by Remark \ref{Rem_tpMinBound} (by Proposition \ref{Prop_dpMintpMin}, since $T$ is dp-minimal, $T$ has no TP-pattern).  We now construct a finite collection of uniform algorithms for defining $\varphi$-types over finite sets.

 Fix $B \subseteq \mathfrak{C}^{\leng(\overline{y})}$ finite and $p \in S_\varphi(B)$.  Let $\mathcal{B}_n$ be defined as above and choose $n$ maximal such that $\mathcal{B}_n \neq \emptyset$.  Since, by Lemma \ref{Lem_tpMinBoundBn}, $\mathcal{B}_{2K-1} = \emptyset$, we have that $n < 2K - 1$.  Fix any $\beta \in \mathcal{B}_n$ $\le_p$-minimal and consider $H(\beta)$ as defined above.  Let $H(\beta) = \{ \gamma_0, ..., \gamma_{m-1} \}$ where $\gamma_{m-1} = \beta$ and, for all $i < j < m$, $\leng(\gamma_i) \le \leng(\gamma_j)$ (so we order the set $H(\beta)$ by length).  Note that $m$ is a function of $n$, say $m = f(n)$.  To be concrete, define $\overline{b}_{i,j} \in B$ for all $i < f(n)$ and $j < n$ so that $\gamma_i = \langle \overline{b}_{i,j} : j < \leng(\gamma_i) \rangle$ for all $i < f(n)$ (for $j$ outside the given range, let $\overline{b}_{i,j}$ be arbitrary).  Let $s \in {}^{f(n) \times n} 2$ be such that
 \[
  p_{ \{ \overline{b}_{i,j} \} } (x) = \left\{ \varphi(x; \overline{b}_{i,j})^{s(i,j)} \right\}
 \]
 for all $i, j$.  So, if $\delta \in {}^B 2$ is the function associated to $p$, then $s(i,j) = \delta(\overline{b}_{i,j})$ for all $i, j$.   Now consider the following algorithm, dependent only on our choices of $n < 2K-1$, $\overline{b}_{i,j}$ for $(i,j) \in f(n) \times n$, and $s \in {}^{f(n) \times n} 2$:

 For all $\overline{b} \in B$, choose $i_0 < f(n)$ minimal so that $\gamma_{i_0}$ $*$-decides $\varphi(x; \overline{b})$.  First of all, such an $i_0$ exists because otherwise, by Lemma \ref{Lem_BasicBn} (v), if $\beta$ does not $*$-decide $\varphi(x; \overline{b})$, then $\beta \concat \langle \overline{b} \rangle \in \mathcal{B}_{n+1}$, contrary to the fact that $\mathcal{B}_{n+1} = \emptyset$.  Second, by Lemma \ref{Lem_CorrectDecisionsonH} and the way we ordered our $\gamma_i$'s, we have that $\gamma_{i_0}$ $*$-decides $\varphi(x; \overline{b})$ correctly.  Therefore, $\varphi(x; \overline{b}) \in p(x)$ if and only if:
 \begin{itemize}
  \item [(i)] $p_{\gamma_{i_0}}(x) \vdash \varphi(x; \overline{b})$ or
  \item [(ii)] $p_{\gamma_{i_0}}(x)$ does not decide $\varphi(x; \overline{b})$ and, for all $\ell < \leng(\gamma_{i_0})$ such that $p_{\gamma_{i_0}, \langle \overline{b}_{i_0, \ell} \rangle}$ decides $\varphi(x; \overline{b})$, we have that
  \[
   p_{\gamma_{i_0}, \langle \overline{b}_{i_0, \ell} \rangle} (x) \vdash \varphi(x; \overline{b}).
  \]
 \end{itemize}
 Notice that all of the above conditions on $\overline{b}$ are expressable as uniform first-order formulas, dependent only on $n < 2K-1$ and $s \in {}^{f(n) \times n} 2$ and defined over the $\overline{b}_{i,j}$'s.  For example, condition (i) (i.e. $p_{\gamma_{i_0}}(x) \vdash \varphi(x; \overline{b})$) holds if and only if
 \[
  \forall x \left( \bigwedge_{j < \leng(\gamma_{i_0})} \varphi(x; \overline{b}_{i_0, j})^{s(i_0, j)} \rightarrow \varphi(x; \overline{b}) \right)
 \]
 holds.  Finally, use this fact to encode this algorithm into a formula $\delta_{n,s}(\overline{y}; \overline{b}_{i,j})_{i < f(n), j < n}$ so that, for all $\overline{b} \in B$,
 \[
  \varphi(x; \overline{b}) \in p(x) \Leftrightarrow \models \delta_{n,s}(\overline{b}; \overline{b}_{i,j})_{i < f(n), j < n}.
 \]
 Since this construction, hence $\delta_{n,s}$, is uniform, $\{ \delta_{n,s} : n < 2K - 1, s \in {}^{f(n) \times n} 2 \}$ is a finite collection of uniform definitions of $\varphi$-types over finite sets.  Thus, by Lemma \ref{Lemma_ManyOneUDTFS}, $\varphi$ has UDTFS.
\end{proof}

This proof only requires that, for some fixed $\varphi(\overline{x}; \overline{y})$, there exists a fixed bound $N < \omega$ such that $\mathcal{B}_N = \emptyset$ for any choice of $B \subseteq \mathfrak{C}^{\leng(\overline{y})}$ finite and $p \in S_\varphi(B)$.  This has more consequences locally.

\begin{thm}\label{Thm_VCDensityOneLocally}
 If $\varphi(\overline{x}; \overline{y})$ is any formula and $N$ is a positive integer such that, for all $B \subseteq \mathfrak{C}^{\leng(\overline{y})}$ with $|B| = N$, we have that $|S_\varphi(B)| \le \frac{N(N+1)}{2}$, then $\varphi$ has UDTFS.
\end{thm}

\begin{proof}
 This follows from the proof of Theorem \ref{Thm_dpMinUDTFS}, except note that, for each $\langle \overline{b}_0, ..., \overline{b}_{N-1} \rangle \in \mathcal{B}_N$, we have that
 \[
  \left| S_\varphi \left( \{ \overline{b}_i : i < N \} \right) \right| > \frac{N(N+1)}{2}
 \]
 by Lemma \ref{Lem_BasicBn} (iv).  Therefore $\mathcal{B}_N = \emptyset$.
\end{proof}

That is to say that all formulas with VC-density less than $2$ have UDTFS.  Noting the relations of dp-minimality to other model-theoretic dividing lines, we get the following corollary to Theorem \ref{Thm_dpMinUDTFS} (see \cite{Adler2008} and \cite{DGL} for definitions):

\begin{cor}\label{Cor_UDTFSForTheories}
 Let $T$ be a complete first-order theory.

 \begin{itemize}
  \item [(i)] If $T$ is VC-minimal, then $T$ has UDTFS.
  \item [(ii)] If $T$ has VC-density one, then $T$ has UDTFS.
 \end{itemize}
\end{cor}

This corollary follows from the fact that VC-minimal implies dp-minimal (see Proposition 9 in \cite{Adler2008}) and VC-density one implies dp-minimal (see Proposition 3.2 in \cite{DGL}).

One should note that UDTFS does not characterize dp-minimality.  There are examples of stable theories that are not dp-minimal (see Theorem 3.5 (iii) in \cite{OUdpmin}, for example).  However, since stable theories have UDTFS, this shows that UDTFS does not imply dp-minimal.

Theorem \ref{Thm_dpMinUDTFS} gives us concrete examples of theories with UDTFS.  For example, the theory of $p$-adic numbers, $\Th(\mathbb{Q}_p)$, is dp-minimal (see Theorem 6.6 in \cite{DGL}).  Therefore, $\Th(\mathbb{Q}_p)$ has UDTFS.  See \cite{DGL}, \cite{Goodrick}, and \cite{Simon} for more examples of theories that are dp-minimal, hence have UDTFS.


\section{Further Questions}

Consider the following picture that demonstrates the relation of UDTFS to other known properties:

\begin{center}
 \begin{tabular}{c c c c c}
   Weakly o-min &               & Stable         &               &            \\
   $\Downarrow$ &               & $\Downarrow$   &               &            \\
   dp-Minimal   & $\Rightarrow$ & UDTFS          & $\Rightarrow$ & Dependent  \\
 \end{tabular}
\end{center}

A question that is still open, due to Laskowski, is the following:

\begin{ques}\label{Conj_DepUDTFS}
 If $T$ is dependent, then does $T$ have UDTFS?
\end{ques}

If Open Question \ref{Conj_DepUDTFS} were true, then this would show that UDTFS is indeed a generalization of definability of types to dependent theories.  This would be an exciting result, but it seems somewhat unlikely.  To date, we have only shown that simple classes of dependent theories have UDTFS.  For example, even if we only assume $\ID(\varphi) = 2$, it is not yet known if $\varphi$ has UDTFS.  One should note, however, that the bound given in Theorem \ref{Thm_VCDensityOneLocally}, namely $|S_\varphi(B)| \le \frac{|B|(|B|+1)}{2}$, is exactly one less than the bound provided by Sauer's Lemma for independence dimension 2.  Since the Sauer bound is tight, there are examples of independence dimension 2 formulas that do not satisfy the conditions of Theorem \ref{Thm_VCDensityOneLocally} (but just barely).

There are a few other properties that hold for dependent formulas but are still open for UDTFS formulas.  For example, if $\varphi(\overline{x}; \overline{y})$ has UDTFS, then does $\varphi^{\opp}(\overline{y}; \overline{x})$ (the formula with opposite partitioning) have UDTFS?  If $T$ is a theory with UDTFS and $T'$ is a reduct of $T$, then does $T'$ have UDTFS?  These statements are true for dependent formulas and theories, so if Open Question \ref{Conj_DepUDTFS} were true, they would certainly be true for UDTFS formulas and theories.

\end{document}